\documentclass[11pt,a4paper]{amsart}
\usepackage{latexsym,enumerate}
\usepackage{amsmath,amsthm,amsopn,amstext,amscd,amsfonts,amssymb}
\usepackage{fullpage}
\usepackage{hyperref}
\usepackage{eucal}
\usepackage{color}

\numberwithin{equation}{section}

\newtheorem{theorem}{Theorem}
\newtheorem{lemma}{Lemma}

\newtheorem{proposition}{Proposition}
\newtheorem{remark}{Remark}

\numberwithin{theorem}{section}
\numberwithin{corollary}{section}
\numberwithin{lemma}{section}
\numberwithin{definition}{section}
\numberwithin{proposition}{section}
\numberwithin{remark}{section}
\newcommand{\eps}{\varepsilon}

\newcommand{\R}{\mathbb R}
\newcommand{\N}{\mathbb N}

\newcommand{\medint}{-\kern  -,375cm\int}

\newcommand{\ie}{\emph{i.e.}}
\newcommand{\eg}{\emph{e.g.}}

\newcommand{\Real}{\mathbb{R}}
\newcommand{\Nat}{\mathbb{N}}

\newcommand{\dist}{\mathop{\mathrm{dist}}\nolimits}
\newcommand{\divergence}{\mathop{\mathrm{div}}\nolimits}
\newcommand{\essinf}{\mathop{\mathrm{ess\;\!inf}}}
\newcommand{\Dom}{\mathsf{D}}
\newcommand{\sii}{L^2}

\usepackage[normalem]{ulem}
\usepackage{color}
\definecolor{DarkBlue}{rgb}{0,0.1,0.7}

\newcommand\soutD{\bgroup\markoverwith
{\textcolor{DarkBlue}{\rule[.5ex]{2pt}{1pt}}}\ULon}
\newcommand{\Hm}[1]{\leavevmode{\marginpar{\tiny%
$\hbox to 0mm{\hspace*{-0.5mm}$\leftarrow$\hss}%
\vcenter{\vrule depth 0.1mm height 0.1mm width \the\marginparwidth}%
\hbox to
0mm{\hss$\rightarrow$\hspace*{-0.5mm}}$\\\relax\raggedright #1}}}

\begin{document}
\title{The equality case in a Poincar\'e-Wirtinger type inequality}
\author{B.~Brandolini, F.~Chiacchio, D.~Krej\v{c}i\v{r}\'ik and C.~Trombetti}

 \address{
Barbara Brandolini \\
Universit\`a degli Studi di Napoli ``Federico II''\\
Dipartimento di Ma\-te\-ma\-ti\-ca e Applicazioni ``R. Caccioppoli''\\
Complesso Monte S. Angelo - Via Cintia\\
80126 Napoli, Italia.
}
\email{brandolini@unina.it}

\address{
Francesco Chiacchio \\
Universit\`a degli Studi di Napoli ``Federico II''\\
Dipartimento di Ma\-te\-ma\-ti\-ca e Applicazioni ``R. Caccioppoli''\\
Complesso Monte S. Angelo - Via Cintia\\
80126 Napoli, Italia.
}
\email{francesco.chiacchio@unina.it}

 \address{David Krej\v{c}i\v{r}\'ik\\
Department of Theoretical Physics, Nuclear Physics Institute \\
Czech Academy of Sciences, \v{R}e\v{z}, Czech Republic.
}
\email{krejcirik@ujf.cas.cz}

 \address{
Cristina Trombetti \\
Universit\`a degli Studi di Napoli ``Federico II''\\
Dipartimento di Ma\-te\-ma\-ti\-ca e Applicazioni ``R. Caccioppoli''\\
Complesso Monte S. Angelo - Via Cintia \\
80126 Napoli, Italia.
}
\email{cristina@unina.it}

\subjclass[2010]{}
\date{}

\begin{abstract}
In this paper, generalizing to the non smooth case already existing results, 
we prove that, for any convex planar set $\Omega$, 
the first non-trivial Neumann eigenvalue $\mu_1(\Omega)$ of 
the Hermite operator 
is greater than or equal to  1.  
Furthermore, and this is our main result, 
under some additional  assumptions on $\Omega$, 
we show that $\mu_1(\Omega)=1$ if and only if $\Omega$ is any strip. 
The study of the equality case requires, among other things, 
an asymptotic analysis of the eigenvalues of 
the Hermite operator in thin domains.
\end{abstract}

\maketitle
\section{Introduction}
%
Let~$\Omega \subset \R^2$ be a convex domain  and let us
denote by~$\gamma$ and~$dm_\gamma$ the standard Gaussian
function and measure in~$\mathbb{R}^{2}$ respectively, that is
\begin{equation*}
  \gamma(x,y):= \exp\left(-\frac{x^2+y^2}{2}\right)
  \qquad \mbox{and} \qquad
  dm_\gamma := \gamma(x,y) \, dx \, dy
  \,.
\end{equation*}
In this paper we consider the following Neumann eigenvalue problem 
for the Hermite operator
\begin{equation}\label{problem}
\left\{
\begin{array}{ll}
  - \divergence (\gamma \nabla u) 
   = \mu \gamma u
  &  \mbox{in} \quad \Omega \,, 
  \\ \\
  \dfrac{\partial u}{\partial \bf{n}} =0
   & \mbox{on} \quad \partial\Omega \,, 
\end{array}
\right.
\end{equation}
where~$\bf{n}$ stands for the outward normal to $\partial \Omega$.
As usual, we understand~\eqref{problem} as a spectral problem for 
the self-adjoint operator~$T$ in the Hilbert space
$\sii_\gamma(\Omega) := \sii(\Omega,dm_\gamma)$ 
associated with the quadratic form 
$
  t[u] := \|\nabla u\|^2
$,
$
  \Dom(t) := H_\gamma^1(\Omega)
$.
Here $\|\cdot\|$ denotes the norm in $\sii_\gamma(\Omega)$ and 
\begin{equation*}
  H_\gamma^1(\Omega) := 
  \{u \in \sii_\gamma(\Omega) \ | \ 
  \nabla u \in \sii_\gamma(\Omega) \}
\end{equation*}
is a weighted Sobolev space equipped with the norm
$\sqrt{\|\cdot\|^2+\|\nabla \cdot\|^2}$. 
Since the embedding 
$H_\gamma^1(\Omega) \hookrightarrow \sii_\gamma(\Omega)$
is compact (see Remark \ref{oss} below), the spectrum of~$T$ is purely discrete.
We arrange the eigenvalues of~$T$ in a non-decreasing sequence
$\{\mu_n(\Omega)\}_{n=0}^{+\infty}$ where each eigenvalue
is repeated according to its multiplicity.
The first eigenfunction of~\eqref{problem} is clearly a constant 
with eigenvalue $\mu_{0}(\Omega)=0$ for any~$\Omega$.
We shall be interested in the first non-trivial eigenvalue 
$\mu _{1}(\Omega)$ of~\eqref{problem},
which admits the following variational characterization
\begin{equation}\label{minimax.intro}
  \mu _{1}(\Omega ) =
  \min
  \left\{ 
  \dfrac{\displaystyle\int_{\Omega }|\nabla u|^{2}\,dm_\gamma}
  {\displaystyle\int_{\Omega } u^{2}\,dm_\gamma} \ : \
  u\in H_\gamma^{1}(\Omega)\setminus \{0\} \,, \
  \int_{\Omega } u \, dm_\gamma=0\right\} 
  .
\end{equation}

A classical Poincar\'e-Wirtinger type inequality 
which goes back to Hermite (see for example \cite[Chapter II, p.~91 ff]{CH1}) 
states that
\begin{equation}\label{poincare}
\mu_1(\R^2)=1
\end{equation}
and therefore
\begin{equation*}
  \int_{\R^2 }\left( u-\int_{\R^2}u \, dm_\gamma\right) ^{2}dm_\gamma
  \leq 
  \int_{\R^2}|\nabla u|^{2} \, dm_\gamma
  \,,
  \qquad \forall u\in H_\gamma^{1}(\R^2)
  \,. 
\end{equation*}
Very recently an inequality analogous to \eqref{poincare} raised up in connection with the proof of the ``gap conjecture'' for bounded sets (see \cite{AC}). In \cite{Andrews-Ni_2012} the authors prove that if $\Omega$ is a bounded, convex set then
\begin{equation}\label{an-ni}
\mu_1(\Omega) \ge \mu_1\left(-\frac{\mathrm{d}(\Omega)}{2},\frac{\mathrm{d}(\Omega)}{2}\right)
\end{equation} 
where $\mathrm{d}(\Omega)$ is the diameter of $\Omega$ and, here and throughout, $\mu_1\left(a,b\right)$ will denote the first nontrivial eigenvalue of the Sturm-Liouville problem
\begin{equation}\label{N_1d}
\left\{ 
\begin{array}{ll}
-\left( \gamma_1v^{\prime }\right) ^{\prime }=\mu \gamma_1v & 
\mbox{in}\ \left( a,b\right) \,, \\ \\
v^{\prime }\left( a\right) =v^{\prime }\left( b\right) =0 \,, &  
\end{array}
\right.  
\end{equation}
with $-\infty \le a <b \le +\infty$ and 
\begin{equation*}
  \gamma_1(x):= \exp\left(-\frac{x^2}{2}\right).
\end{equation*}
Again, we understand~\eqref{N_1d} 
as a spectral problem for a self-adjoint operator 
with compact resolvent in $\sii_{\gamma_1}((a,b))$.
It is well-known that
\begin{equation}\label{l(R)}
\mu _{1}(a,b)\geq 1
\quad \text{ with } \quad \mu _{1}(a,b)=1
\text{ \ if and only if \ }(a,b)=\mathbb{R}.  
\end{equation}

As first result of this paper we extend the validity of~\eqref{an-ni}
to any convex, possibly unbounded, planar domain (see \cite{BCHT} for the smooth case).
\begin{theorem}\label{main copy(1)} 
Let $\Omega \subset \mathbb{R}^{2}$ be any convex domain. Then
\begin{equation}\label{mi}
\mu _{1}(\Omega )\geq \mu_1(\R)=1 \,.
\end{equation}
\end{theorem}

The result is sharp in the sense that 
the equality in~\eqref{mi} is achieved 
for~$\Omega$ being any two-dimensional strip.
It is natural to ask if the strips are the unique domains
for which the equality in~\eqref{mi} is achieved. 
 
We provide a partial answer to the uniqueness question
via the following theorem, which is the main result of this paper. 
\begin{theorem}\label{main}
Let $\Omega $ be a convex  subset of 
$
  S_{y_{1},y_{2}} := 
  \left\{ (x,y)\in \mathbb{R}^{2}:y_{1}<y<y_{2}\right\}
$ 
for some $y_1,y_2\in \R$, $y_1 < y_2$.
If $\mu _{1}(\Omega )=1$, then $\Omega $ is a strip.
\end{theorem}

Inequality \eqref{an-ni} is a Payne-Weinberger type inequality 
for the Hermite operator. We recall that the classical Payne-Weinberger inequality 
states that the first nontrivial eigenvalue of the Neumann Laplacian 
in a bounded convex set~$\Omega$, $\mu_1^{\Delta}(\Omega)$, 
satisfies the following bound
\begin{equation}\label{PW}
\mu_1^{{\Delta}}(\Omega)\ge \frac{\pi^2}{\mathrm{d}(\Omega)^2},
\end{equation}
where $\pi^2/\mathrm{d}(\Omega)^{2}$  is the first
nontrivial Neumann eigenvalue of the one-dimensional Laplacian 
in $\left( -\mathrm{d}(\Omega)/2,\mathrm{d}(\Omega)/2\right)$ (see \cite{PW}). 
The above estimate is the best bound that can be given in terms of 
the diameter alone in the sense that 
$\mu_1^{\Delta}(\Omega){\mathrm{d}(\Omega)^2}$ tends to $\pi^2$ for a
parallepiped all but one of whose dimensions shrink to zero (see \cite{HW,V}).

Estimate \eqref{mi} is sharp, not only asymptotically, 
since the equality sign is achieved when $\Omega$ is any strip~$S$. 
Indeed, it is straight-forward to verify that $\mu_1(S)=\mu_1(\R)=1$ for any strip~$S$.
Hence the question faced in Theorem \ref{main} appears quite natural.

The paper is organized as follows. 
Section~\ref{Sec.ext} contains the proof of Theorem~\ref{main copy(1)}, 
while Section~\ref{Sec.main} is devoted to the proof of Theorem~\ref{main}. 
The latter consists in various steps.  
We firstly deduce from~\eqref{mi} that any optimal set must be unbounded; 
then we show that it is possible to split an optimal set~$\Omega$ 
getting two sets that are  still optimal and have Gaussian area $m_\gamma(\Omega)/2$. 
Repeating  this procedure we obtain a sequence of thinner and thinner, 
optimal sets~$\Omega_k$ and
we finally prove that there exists $a\in \overline \R$ 
such that $\mu_1(\Omega_k)$ converges as $k\to+\infty$ to $\mu_1(a,+\infty)$, 
which is strictly greater than 1 unless $a=-\infty$. 
This circumstance implies that $\Omega$ contains a straight-line, 
and hence~$\Omega$ is a strip. 

The convergence of $\mu_1(\Omega_k)$ to $\mu_1(a,+\infty)$
follows by a more general result established in Section~\ref{Sec.as},
where we actually prove a convergence of \emph{all} eigenvalues of~$T$
in thin domains to eigenvalues of a one-dimensional problem
(see ~Theorem~\ref{Thm.convergence}).
We also establish certain convergence of eigenfunctions.
We believe that the convergence results are of independent interest,
since our method of proof differs from known techniques
in the case of the Neumann Laplacian in thin domains
\cite{Arrieta_1995,Arrieta-Carvalho_2004,Jerison2,Schatzman_1996}.

For optimisation results related to the present work,
we refer the interested reader to \cite{BCT_p,ENT, FNT,BCM,CdB}.

\section{Proof of Theorem~\ref{main copy(1)}}\label{Sec.ext}
%
Repeating step by step the arguments contained in~\cite{BCHT}, 
Theorem~\ref{main copy(1)} is a consequence of the following extension result,
which we believe is interesting on its own.
\begin{theorem}\label{ext}
Let $\Omega\subset \R^2$ be a convex domain
and let $u \in H_\gamma^1(\Omega)$. 
Then there exists a function $\tilde u \in H_\gamma^1(\R^2)$ 
which is an extension of $u$ to $\R^2$, that is
$$
\tilde u|_{\Omega}=u
$$
and
\begin{equation}\label{extension}
  \|\tilde u\|_{H_\gamma^1(\R^2)}\le C \, \|u\|_{H_\gamma^1(\Omega)}
  \,,
\end{equation}
where $C=C(\Omega)$.
\end{theorem}

\begin{proof}
We preliminarily observe that, if $\Omega$ is bounded, 
the theorem can be immediately obtained from the classical result 
for the unweighted case 
(see for instance \cite[Thm.~4.4.1]{Evans-Gariepy}). 
So, from now on, we assume that $\Omega$ is unbounded and we adapt 
the arguments in \cite{Evans-Gariepy} to treat our case. 
We distinguish two cases: $0 \in \Omega$ and $0 \notin \Omega$. 

\medskip
\fbox{\emph{Case 1:} $0 \in \Omega$.} 
The convexity of $\Omega$ ensures there exists a constant $L>0$ such that, 
for every $(x_0,y_0)\in \partial\Omega$, up to a rotation, 
there exist $r>0$ and an $L$-Lipschitz continuous function 
$\beta:\R \to [0,+\infty)$ such that, 
if we set $Q(x_0,y_0,r):=\{(x,y)\in \R^2:\> |x-x_0|<r,\> |y-y_0|<r\}$, 
it holds
$$
\Omega \cap Q(x_0,y_0,r)=\{(x,y)\in \Omega: \> |x-x_0|<r,\> y_0-r<y<\beta(x)\}, \quad \max_{|x-x_0|<r}|\beta(x)-y_0|<\frac{r}{2}.
$$
In other words, 
\begin{equation}\label{beta'}
|\beta'(x)|\le L \qquad \mbox{for a.e.} \ \, x\in (x_0-r,x_0+r),
\end{equation}
with $L$ independent from $x_0,y_0,r$. 

Fix $(x_0,y_0)\in \partial\Omega$ and set 
$\Omega^i:=Q(x_0,y_0,r)\cap \Omega$ and 
$\Omega^e:=Q(x_0,y_0,r)\setminus \overline{\Omega}$. 
Let $u \in C^1(\overline{\Omega})$ and suppose for the moment that the support of $u$ is contained in 
$ Q(x_0,y_0,r)\cap \overline{\Omega}$. 
Set
$$
u^e(x,y) := u(x,2\beta(x)-y) 
\qquad \mathrm{if} \> (x,y)\in \overline{\Omega^e}\,.
$$
We get
\begin{eqnarray*}
\lefteqn{\int_{\Omega^e}{u^e(x,y)^2} \, \exp\left(-\frac{x^2+y^2}{2}\right)dxdy}
\\ 
&&=\int_{\Omega^e}u(x,2\beta(x)-y)^2 \, \exp\left(-\frac{x^2+y^2}{2}\right) dxdy
\\
&&= \int_{\Omega^i} u(s,t)^2 \, \exp\left(-\frac{s^2+(2\beta(s)-t)^2}{2}\right)dsdt
\\
&&=\int_{\Omega^i} u(s,t)^2 \, \exp\left(-\frac{s^2+(2\beta(s)-t)^2}{2}
+\frac{s^2+t^2}{2}\right) \, \exp\left(-\frac{s^2+t^2}{2}\right)dsdt.
\end{eqnarray*}
By elementary geometric considerations, taking into account the assumption $0\in \Omega$, it is easy to verify that
\begin{equation}
\exp\left(-\frac{s^2+(2\beta(s)-t)^2}{2}+\frac{s^2+t^2}{2}\right) \leq 1
  \,, \qquad
  \forall \>(s,t)\in \Omega^i.  \label{exp}
\end{equation}
Thus
\begin{equation}\label{normau}
\displaystyle\int_{\Omega^e}{u^e(x,y)}^{2}
\, \exp\left(-\frac{x^2+y^2}{2}\right)dxdy 
\le \int_{\Omega^i} u(s,t)^2 \, \exp\left(-\frac{s^2+t^2}{2}\right)dsdt.
\end{equation}
On the other hand, by \eqref{beta'} and \eqref{exp} it holds
\begin{eqnarray}\label{normaDu}
\lefteqn{\int_{\Omega^e}|\nabla u^e(x,y)|^2
\,\exp\left(-\frac{x^2+y^2}{2}\right)dxdy} 
\notag
\\
&&\le \int_{\Omega^i} 
\left[\left(\partial_s u(s,t)+2 \partial_t u(s,t)
\beta'(s)\right)^2+\big(\partial_t u(s,t)\big)^2\right]
\exp\left(-\frac{s^2+(2\beta(s)-t)^2}{2}\right)dsdt 
\notag
\\
&&\le C(L)\int_{\Omega^i} |\nabla u(s,t)|^2
\,\exp\left(-\frac{s^2+t^2}{2}\right)dsdt. 
\end{eqnarray}
Define
$$
\tilde u:=\left\{
\begin{array}{ll}
u & \mathrm{on}\quad \overline{\Omega^i} \,,
\\
u^e & \mathrm{on}\quad \overline{\Omega^e} \,,  
\\
0 & \mathrm{on} \quad \R^2 \setminus(\Omega^i\cup\Omega^e) \,.
\end{array}
\right.
$$
If~$\Omega $ contains the origin and the support of $u$ is contained in
$Q(x_0,y_0,r)\cap \overline{\Omega}$, 
then~\eqref{normau} and~\eqref{normaDu} imply~\eqref{extension} with $C=C(L)$.

Now assume that $u \in C^1(\overline{\Omega})$ 
and drop the restriction on its support. 
Clearly $\partial \Omega$ is not compact, 
but we can cover $\partial \Omega$  with a countable family of squares 
$\{Q_{2k-1}\}=\{Q(x_0^k,y_0^k,r'_k)\}_{k=1}^\infty$, 
with $(x_0^k,y_0^k)\in\partial\Omega$ and $r'_k>0$ 
such that \eqref{beta'} holds true for each~$k$. 
Analogously, we can cover the set 
$\Omega\setminus \bigcup_{k=1}^{+\infty} Q_{2k-1}$ 
with a countable family of squares 
$\{Q_{2k}\}=\{Q(x_1^k,y_1^k,r''_k)\}_{k=1}^\infty$ 
with $(x_1^k,y_1^k)\in\Omega$ and $r''_k>0$. 
For the countable cover $\{Q_l\}_{l=1}^\infty$ of $\overline{\Omega}$ 
there is a partition of unity $\varphi_l$  subordinated to $Q_l$ with $\varphi_l$ 
smooth for each $l$ (see for instance~\cite[Thm. 3.14]{A}). 
Define $\widetilde{\varphi_l u}$ as above when~$l$ is odd and set
$$
\tilde u := \sum_{l \>\mathrm{odd}}\widetilde{\varphi_l u}
+\sum_{l \>\mathrm{even}} \varphi_l u.
$$
Clearly $\tilde u$ satisfies \eqref{extension}.

Finally if $u \in H_\gamma^1(\Omega)$, 
the claim follows by approximation arguments 
and the proof of Case~1 is accomplished. 

\medskip
\fbox{\emph{Case 2:} $0 \not\in \Omega$.}
Suppose now that $0\notin \Omega $ and denote 
$d_{0} := \dist(0,\partial \Omega )$. 
Let us fix a vector $(\delta_1, \delta_2)$,  
with $\sqrt{\delta_1^2+\delta_2^2}>d_0$, in such a way that the translation
$
  \Phi: \R^2 \to \R^2:
  \{(x,y) \mapsto (x-\delta_1,y-\delta_2) \}
$ 
maps $\Omega $ onto a set $\Phi(\Omega )$ containing the origin. 
Defining
\begin{equation*}
v(x,y) := u(x+\delta_1,y+\delta_2)\exp \left( -
\frac{x\delta_1 }{2}-\frac{\delta_1 ^{2}}{4}\right)\exp \left( -
\frac{y\delta_2 }{2}-\frac{\delta_2 ^{2}}{4}\right)  
, 
\end{equation*}
for every $(x,y)\in \Phi(\Omega)$,
we have
\begin{equation*}
\int_{\Omega }u^{2} \, dm_\gamma = \int_{\Phi(\Omega )} v^{2} \, dm_\gamma
\,.  
\end{equation*}
Since by construction $\Phi(\Omega )$ contains the origin, 
there exists a function $\tilde{v}\in H_\gamma^{1}(\mathbb{R}^{2})$ 
such that 
$\tilde{v}\left\vert _{\Phi(\Omega )}\right. =v$ 
and
\begin{equation*}
\|\tilde{v}\|_{H_\gamma^{1}(\mathbb{R}^{2})}
\leq C(L) \, \|v\|_{H_\gamma^{1}(\Phi(\Omega ))}
\,.
\end{equation*}
Letting
\begin{equation*}
\tilde{u}(x,y) :=\tilde{v}(x-\delta_1
,y-\delta_2)\exp \left( \frac{x\delta_1 }{2}-\frac{\delta_1 ^{2}}{4}\right)
\exp \left( \frac{y\delta_2 }{2}-\frac{\delta_2 ^{2}}{4}\right) 
,
\end{equation*}
we finally get that $\tilde u|_\Omega=u$ and
\begin{equation*}
\|\tilde{u}\|_{H_\gamma^{1}(\mathbb{R}^{2})}
\leq C(L,d_{0}) \, \|u\|_{H_\gamma^{1}(\Omega)}.
\end{equation*}
This completes the proof of the theorem.
\end{proof}

\begin{remark}\label{oss}
Using the fact that $H^{1}_\gamma(\mathbb{R}^{2})$ is compactly
embedded into $L^{2}_\gamma(\mathbb{R}^{2})$ 
(see for example~\cite{D}) and the above extension theorem 
one can easily deduce the compact embedding of 
$H^{1}_\gamma(\Omega)$ into $L^{2}_\gamma(\Omega)$ (see also \cite{FP,BCHT}).  
Therefore, by the classical spectral theory on compact self-adjoint
operators, $\mu _{1}(\Omega )$ satisfies the variational 
characterization~\eqref{minimax.intro}.
\end{remark}

\section{Proof of Theorem~\ref{main}}\label{Sec.main}

The main ingredient in our proof of Theorem~\ref{main} 
is the following lemma, which  tells us that
cutting the optimiser of~\eqref{mi} in two convex, unbounded sets with equal Gaussian area, 
we again get two optimisers.
\begin{lemma}
Let $\Omega$ be a convex subset of $S_{y_1,y_2}$
 with $\mu_1(\Omega)=1$. 
Let $\bar y\in (y_1,y_2)$ be such that the straight-line 
$\{y=\bar y\}$ divides $\Omega$ into two convex subsets with equal Gaussian area $m_\gamma(\Omega)$. 
Then
$$
\mu_1\big(\Omega \cap \{y<\bar y\}\big)
=\mu_1\big(\Omega \cap \{y>\bar y\}\big)=1.
$$
\end{lemma}
\begin{proof}
Let $u$ be an eigenfunction of
(\ref{problem}) corresponding to $\mu _{1}(\Omega ).$ 
By~\eqref{minimax.intro}, we know that $\int_{\Omega} u \, dm_\gamma=0$ and 
\begin{equation*}
1=\dfrac{\int_{\Omega }|\nabla u|^{2}\,dm_\gamma}
{\int_{\Omega }u^{2}\,dm_\gamma}
\,.
\end{equation*}
For each $\alpha \in [0,2\pi]$ there is a unique straight-line $r_\alpha$ 
orthogonal to $(\cos \alpha,\sin \alpha)$ 
such that it divides $\Omega$ into two convex sets 
$\Omega_\alpha',\Omega_\alpha''$ with equal Gaussian measure. 
Let $I(\alpha):=\int_{\Omega_\alpha'}u \, dm_\gamma$.  
Since $I(\alpha)=-I(\alpha+\pi)$, 
by continuity there is $\bar \alpha $ such that $I(\bar \alpha)=0$. 
Now we claim that $r_{\bar \alpha}$ is parallel to the $x$-axis. 
Note firstly that $\Omega_{\bar \alpha}'$ and $\Omega_{\bar \alpha}''$ 
are obviously convex and by \eqref{an-ni}, \eqref{l(R)} and \eqref{mi} we have
\begin{equation}
\mu _{1}(\Omega_{\bar \alpha}')\geq 1 \,,
\qquad \mu_1(\Omega_{\bar \alpha}'')\ge 1 \,.  
\label{mu+->=1}
\end{equation}
Moreover, it is immediate to verify that
\begin{equation*}
1=\mu_1(\Omega)
=\dfrac{\int_{\Omega_{\bar \alpha}'}|\nabla u|^{2}\,dm_\gamma
+\int_{\Omega_{\bar \alpha}''}|\nabla u|^{2}\,dm_\gamma}
{\int_{\Omega_{\bar \alpha}'}u^{2}\,dm_\gamma
+\int_{\Omega_{\bar \alpha}''}u^{2}\,dm_\gamma}
\geq \min \left\{ \dfrac{\int_{\Omega_{\bar \alpha}'}|\nabla u|^{2}\,dm_\gamma}
{\int_{\Omega_{\bar \alpha}'}u^{2}\,dm_\gamma},
\dfrac{\int_{\Omega_{\bar \alpha}''}|\nabla u|^{2}\,dm_\gamma}
{\int_{\Omega_{\bar \alpha}''}u^{2}\,dm_\gamma}\right\},
\end{equation*}
with equality holding if and only if
$$
\dfrac{\int_{\Omega_{\bar \alpha}'}|\nabla u|^{2}\,dm_\gamma}
{\int_{\Omega_{\bar \alpha}'}u^{2}\,dm_\gamma}
=\dfrac{\int_{\Omega_{\bar \alpha}''}|\nabla u|^{2} \, dm_\gamma}{
\int_{\Omega_{\bar \alpha}''}u^{2}\,dm_\gamma}.
$$
Without loss of generality we can assume that
\begin{equation*}
\min \left\{ \dfrac{\int_{\Omega_{\bar \alpha}'}|\nabla u|^{2}\,dm_\gamma}
{\int_{\Omega_{\bar \alpha}'}u^{2}\,dm_\gamma},
\dfrac{\int_{\Omega_{\bar \alpha}''}|\nabla u|^{2}\,dm_\gamma}
{\int_{\Omega_{\bar \alpha}''}u^{2}\,dm_\gamma}\right\} 
=\dfrac{\int_{\Omega_{\bar \alpha}'}|\nabla u|^{2}\,dm_\gamma}
{\int_{\Omega_{\bar \alpha}'}u^{2}\,dm_\gamma}.
\end{equation*}
Finally, (\ref{mu+->=1}) ensures that
\begin{equation}
1=\mu _{1}(\Omega )=\mu _{1}(\Omega_{\bar \alpha}')
=\mu _{1}(\Omega_{\bar \alpha}'') \,.  \label{1=m}
\end{equation}

Now we want to show that both $\Omega_{\bar \alpha}'$ and $\Omega_{\bar \alpha}''$ are unbounded, and hence $r_{\bar \alpha}$ is parallel to the $x$-axis.
Suppose by contradiction that, for instance, $\Omega_{\bar \alpha}'$ is bounded. In
such a case \eqref{an-ni} yields 
\begin{equation*}
\mu _{1}(\Omega_{\bar \alpha}')\geq \mu _{1}\left( -\frac{\mathrm{d}(\Omega_{\bar \alpha}')}{2},
\frac{\mathrm{d}(\Omega_{\bar \alpha}')}{2}\right) .
\end{equation*}
Taking into account~(\ref{1=m}) and (\ref{l(R)}), 
we get that
\begin{equation*}
\mu_{1}\left( -\frac{\mathrm{d}(\Omega_{\bar \alpha}')}{2},\frac{\mathrm{d}
(\Omega_{\bar \alpha}')}{2}\right) =1
\end{equation*}
 that is $\mathrm{d}(\Omega_{\bar \alpha}')=+\infty $, which is a
contradiction. 
\end{proof}
\begin{proof}[Proof of Theorem \ref{main}]

By contradiction, let us assume that $\Omega \subset S_{y_1,y_2}$ 
is a convex domain different from a strip
and $\mu_1(\Omega)=1$. Let us denote
$$
\Omega=\left\{(x,y)\in \R^2: \, y_1<y<y_2,\, p(y)<x\right\},
$$
where $p$ is a convex, non-trivial function. 
From~\eqref{an-ni} and~\eqref{l(R)}
it follows that~$\Omega$ is necessarily unbounded. 
By employing a separation of variables, 
we also deduce from~\eqref{an-ni}  and~\eqref{l(R)}
that~$\Omega$ cannot be a semi-strip. 
Finally, we may assume that 
$
  \inf\{x: \exists\, y\in[y_1,y_2], \, (x,y)\in\Omega\}
$
is finite
(otherwise, we would have the finite supremum,
which can be transferred to our situation by a reflection
of the coordinate system). 

Repeating the procedure described in the above lemma, 
since at any step we are dividing into two convex subsets with equal Gaussian area, 
we can obtain a sequence of unbounded convex domains 
\begin{equation}\label{Omega.cut}
  \Omega_{\epsilon_k} 
  := \left\{ (x,y)\in \mathbb{R}^{2}:\ y_0<y<d_k,p(y)<x \right\}
  = \left\{ (x,y)\in \Omega:\ y_0<y<d_k \right\}
\end{equation}
such that
$$
  \mu _{1}(\Omega_{\epsilon_k}) =1 
  \,, \qquad
  \epsilon_{k} := d_{k}-y_0 \xrightarrow[k \to +\infty]{} 0
  \,.
$$
Here the point~$y_0$ is chosen in such a way that  
$
  p^{\prime}(y_0) \not= 0
$,
which is always possible because the situation of semi-strips
has been excluded.
Without loss of generality 
(reflecting again the coordinate system if necessary), 
we may in fact assume 
\begin{equation}\label{Ass.phi}
  p^{\prime}(y_0) > 0
  \,,
\end{equation}
so that~$\phi$ is increasing on $[y_0,d_k]$
whenever~$k$ is sufficiently large.
Applying now a more general convergence result for eigenvalues
in thin Neumann domains that we shall establish 
in the following section (Theorem~\ref{Thm.convergence}), 
we have 
\begin{lemma}\label{Lem.convergence}
$
\displaystyle
  \lim_{k\to\infty}\mu_1(\Omega_{\epsilon_k})
  =\mu_1(p^{-1}(y_0),+\infty)
$.
\end{lemma}
\noindent
Since $\mu _{1}(\Omega_{\epsilon_k})$ equals~$1$ for every~$k$,
we conclude that
$$
  \mu_1\big(p^{-1}(y_0),+\infty\big) = 1
  \,.
$$
However, from \eqref{l(R)}, 
we then deduce that $p^{-1}(y_0)=-\infty$,
which contradicts our assumptions
from the beginning of the proof.
In other words, $\Omega$ contains a straight-line 
and the theorem immediately follows. 
\end{proof}

It thus remains to establish Lemma~\ref{Lem.convergence}.

\section{Eigenvalue asymptotics in thin strips}\label{Sec.as}
%
In this section we establish Lemma~\ref{Lem.convergence}
as a consequence of a general result about convergence
of \emph{all} eigenvalues of~$T$
in thin domains of the type~\eqref{Omega.cut}.

\subsection{The geometric setting}
%
Let $f:[0,+\infty) \to [0,+\infty)$ be a concave 
non-decreasing continuous non-trivial function such that $f(0) = 0$
(the case $f(0)>0$ is actually much easier to deal with). 
Given a positive number $\eps < \sup f$,
we put
$$
  f_\eps(x) := \min\{\eps,f(x)\}
$$
and define an unbounded domain
$$
  \Omega_\eps := \{(x,y)\in\Real^2 :
  0 < x \,, \ 0<y<f_\eps(x)\}
  \,.
$$
Clearly, \eqref{Omega.cut} can be cast into this form after 
identifying $f=p^{-1}$ and a translation. 
However, keeping in mind that the problem~\eqref{problem} 
is not translation-invariant, we accordingly change 
the definition of the Gaussian weight
throughout this section
$$
  \gamma(x,y) := \exp\left(-\frac{(x_0+x)^2+(y_0+y)^2}{2}\right)
  \,.
$$
Here~$y_0$ is primarily thought as the point from~\eqref{Omega.cut}
and~$x_0$ is then such that $(x_0,y_0) \in \Omega_{\epsilon_k}$. 
For the results established in this section, however, 
$x_0$ and $y_0$ can be thought as arbitrary real numbers.
For our method to work, it is only important to assume~\eqref{Ass.phi},
which accordingly transfers to
\begin{equation}\label{Ass}
  f'(0) < +\infty \,.
\end{equation}
%

\subsection{The analytic setting and main result}
%
Keeping the translation we have made in mind,
instead of~\eqref{problem} we equivalently consider the eigenvalue problem 
\begin{equation}\label{bv}
\left\{
\begin{array}{ll}
  - \divergence (\gamma \nabla u) 
  = \mu \gamma u
  & \mbox{in} \quad \Omega_\eps \,, 
  \\ \\
  \dfrac{\partial u}{\partial \bf{n}} 
  =0
  & \mbox{on} \quad \partial\Omega_\eps \,.
\end{array}
\right.
\end{equation}
We understand~\eqref{bv} as a spectral problem for 
the self-adjoint operator~$T_\eps$ in the Hilbert space
$\sii_\gamma(\Omega_\eps)$ 
associated with the quadratic form 
$
  t_\eps[u] := \|\nabla u\|_\eps^2
$,
$
  \Dom(t_\eps) := H_\gamma^1(\Omega_\eps)
$.
Here $\|\cdot\|_\eps$ denotes the norm in $\sii_\gamma(\Omega_\eps)$.
We arrange the eigenvalues of~$T_\eps$ in a non-decreasing sequence
$\{\mu_n(\Omega_\eps)\}_{n \in \Nat}$ where each eigenvalue
is repeated according to its multiplicity.
In this paper we adopt the convention $0 \in \Nat$.
We are interested in the behaviour of the spectrum as $\eps \to 0$,
particularly $\mu_1(\Omega_\eps)$ because of Lemma~\ref{Lem.convergence}.   
 
It is expectable that the eigenvalues will be determined 
in the limit $\eps \to 0$ by the one-dimensional problem
\begin{equation}\label{bv.1D}
\left\{
\begin{array}{ll}
  - (\gamma_0 \, u')' 
  = \nu \gamma_0 u
  & \mbox{in} \quad (0,+\infty) \,, 
  \\ \\
  u' (0)=0,
    & 
\end{array}
\right.
\end{equation}
where
$$
  \gamma_0(x) := \gamma(x,0)
  = \exp\left(-\frac{(x_0+x)^2+y_0^2}{2}\right)
  \,.
$$
Again, we understand~\eqref{bv.1D} as a spectral problem for 
the self-adjoint operator~$T_0$ in the Hilbert space
$\sii_{\gamma_0}((0,+\infty))$ 
associated with the quadratic form 
$
  t_0[u] := \|\nabla u\|_0^2
$,
$
  \Dom(t_0) := H_{\gamma_0}^1((0,+\infty))
$, 
where $\|\cdot\|_0$ denotes the norm in $\sii_{\gamma_0}((0,+\infty))$.
As above, we arrange the eigenvalues of~$T_0$ in a non-decreasing sequence
$\{\nu_n\}_{n \in \Nat}$ where each eigenvalue
is repeated according to its multiplicity. By construction,  for each $n\in\N$, $\nu_n$ coincides with the eigenvalue $\mu_n(x_0,+\infty)$ defined in \eqref{N_1d}.

In this section we prove the following convergence result.
\begin{theorem}\label{Thm.convergence}
Let $f:[0,+\infty) \to [0,+\infty)$ be a concave 
non-decreasing continuous non-trivial function such that $f(0) = 0$.
Assume in addition~\eqref{Ass}. 
Then 
$$
  \forall n\in\Nat \,, \qquad
  \mu_n(\Omega_\eps) \xrightarrow[\eps \to 0]{} \nu_n
  \,.
$$ 
\end{theorem}

We shall also establish certain convergence
of eigenfunctions of~$T_\eps$ to eigenfunctions of~$T_0$.

Clearly, Lemma~\ref{Lem.convergence} is the case $n=1$
of this general theorem.

The rest of this section is devoted to 
a proof of Theorem~\ref{Thm.convergence}.

\subsection{From the moving to a fixed domain}
%
Our main strategy is to map~$\Omega_\eps$ into a fixed strip~$\Omega$.
We introduce a refined mapping in order to effectively 
deal with the singular situation $f(0)=0$.

Let
$$
  a_\eps := \inf f_\eps^{-1}(\{\eps\})
  \,.
$$
By the definition of~$f_\eps$ and since~$f$ is non-decreasing, 
$a_\eps \to 0$ as $\eps \to 0$ and
$f_\eps(x)=\eps$ for all $x>a_\eps$.  
If $f(0)>0$, then there exists $\eps_0>0$ such that
$a_\eps=0$ for all $\eps \leq \eps_0$.  
On the other hand, if $f(0)=0$, 
then $a_\eps>0$ for all $\eps>0$.
The troublesome situation is the latter, 
to which we have restricted from the beginning.
In this case, we introduce an auxiliary function
$$
  g_\eps(s) := 
  \begin{cases}
    a_\eps s + a_\eps 
    &\mbox{if} \quad s \in [-1,0) \,,
    \\
    s+a_\eps  
    &\mbox{if} \quad s \in [0,+\infty) \,.
  \end{cases}
$$
Since we are interested in the limit $\eps \to 0$, we may henceforth assume 
\begin{equation}\label{Ass.small}
  \eps \leq 1 
  \qquad \mbox{and} \qquad
  a_\eps \leq 1
  \,.
\end{equation}

Define $\eps$-independent sets 
$$
  \Omega_-:=(-1,0)\times(0,1) 
  \,, \qquad
  \Omega_+:=(0,+\infty)\times(0,1)
  \,, \qquad
  \Omega:=(-1,+\infty)\times(0,1)
  \,.
$$
The mapping
\begin{equation}\label{layer}
  \mathcal{L}_\eps : \Omega \to \Omega_\eps :
  \left\{
  (s,t) \mapsto
  \mathcal{L}_\eps(s,t) := 
  \big(
  g_\eps(s), f_\eps(g_\eps(s)) \, t
  \big)
  \right\}
\end{equation}
represents a $C^{0,1}$-diffeomorphism 
between~$\Omega$ and~$\Omega_\eps$
($f$~is differentiable almost everywhere,
as it is supposed to be concave). 
In this way, {we obtain a convenient parameterisation of~$\Omega_\eps$
via the coordinates $(s,t) \in \Omega$}
whose Jacobian is 
\begin{equation}\label{Jacobian}
  j_\eps(s,t)
  = g_\eps'(s) f_\eps(g_\eps(s))
  \,.
\end{equation}
Note that the Jacobian is independent of~$t$
and singular at $s=-1$.
Now we reconsider~\eqref{bv} in~$\Omega$.
With the notation 
$$
  \gamma_\eps(s,t) := (\gamma \circ \mathcal{L}_\eps)(s,t)
  = \exp\left(-\frac{[x_0+g_\eps(s)]^2
  +[y_0+f_\eps(g_\eps(s))t]^2}{2}\right)
  \,,
$$ 
introduce the unitary transform
$$
  U_\eps: 
  \sii_\gamma(\Omega_\eps) \to 
  \sii_{\gamma_\eps j_\eps/\eps}(\Omega) :
  \left\{
  u \mapsto \sqrt{\eps} \, u \circ \mathcal{L}_\eps
  \right\}
  \,.
$$
Here, in addition to the change of variables~\eqref{layer},
we also make an irrelevant scaling transform 
(so that the renormalised Jacobian~$j_\eps/\eps$ is~$1$ in~$\Omega_+$).
The operators $H_\eps := U_\eps T_\eps U_\eps^{-1}$
and~$T_\eps$ are isospectral.
By definition, $H_\eps$~is associated with the quadratic form
$h_\eps[\psi] := t_\eps[U_\eps^{-1}\psi]$, 
$\Dom(h_\eps) := U_\eps \Dom(t_\eps)$.
\begin{proposition}\label{Prop.form}
Assume~\eqref{Ass}.
Then
\begin{align}
  h_\eps[\psi] &= \int_\Omega 
  \left[
  \left(
  \frac{\partial_s \psi}{g_\eps'} 
  - \frac{f_\eps'\circ g_\eps}{f_\eps \circ g_\eps} \, t \, \partial_t \psi
  \right)^2 
  + 
  \frac{(\partial_t \psi)^2}{(f_\eps\circ g_\eps)^2} 
  \right]
  \gamma_\eps \, g_\eps' \, \frac{f_\eps \circ g_\eps}{\eps} \, ds \, dt
  \,, \qquad
  \label{form1}
  \\
  \Dom(h_\eps) &\subset H_{\gamma_\eps j_\eps/\eps}^1(\Omega)
  \,.
  \label{form2}
\end{align}
\end{proposition}
\noindent
Here we have started to simplify the notation 
by suppressing arguments of the functions.
\begin{proof}
The space
$
  \mathcal{D}_\eps := C_0^1(\Real^2) \upharpoonright\Omega_\eps
$
is a core of~$t_\eps$. 
The transformed space 
$
  \mathcal{D} := U_\eps \mathcal{D}_\eps
$
is a subset of $C_0^0(\Real^2) \upharpoonright \Omega$ consisting of 
Lipschitz continuous functions on~$\Omega$ which belong to 
$C^1(\overline{\Omega_-}) \oplus C^1(\overline{\Omega_+})$
(we do not have $C^1$ globally, 
because~$g_\eps$ and~$f_\eps$ are not smooth).
For any $\psi \in \mathcal{D}$, it is easy to check~\eqref{form1};
this formula extends to all~$\psi$ from the domain
$$
  \Dom(h_\eps) = \overline{\mathcal{D}}^{\|\cdot\|_{h_\eps}}
  \,, \qquad
  \|\cdot\|_{h_\eps} := \sqrt{h_\eps[\cdot] + \|\cdot\|^2}
  \,,
$$ 
where $\|\cdot\|$ denotes the norm of 
$\sii_{\gamma_\eps j_\eps/\eps}(\Omega)$.
Let $\psi \in \mathcal{D}$.
Using elementary estimates, we easily check
\begin{equation}\label{bound.form}
  h_\eps^-[\psi] \leq h_\eps[\psi] 
\end{equation}
where
\begin{align*}
  h_\eps^-[\psi] &:= 
  \delta \int_{\Omega}
  \left(
  \frac{\partial_s \psi}{g_\eps'} 
  \right)^2 
  \gamma_\eps \, g_\eps' \, \frac{f_\eps \circ g_\eps}{\eps} \, ds \, dt
  + \left(
  1 - \frac{\delta}{1-\delta} \, \|f_\eps'\|_\infty^2
  \right) 
  \int_{\Omega} \frac{(\partial_t \psi)^2}{(f_\eps \circ g_\eps)^2}
  \, \gamma_\eps \, g_\eps' \, \frac{f_\eps \circ g_\eps}{\eps} \, ds \, dt
\end{align*}
with any $\delta \in (0,1)$.
Note that~$f_\eps'$ is bounded under the assumption~\eqref{Ass} 
and the concavity.
For any $\eps > 0$,
we can choose~$\delta$ so small that $h_\eps^-[\psi]$
is composed of a sum of two non-negative terms
($\delta$~can be made independent of~$\eps$
if we restrict the latter to a fixed bounded interval,
say $(0,1]$, see \eqref{Ass.small}, because 
$
  \|f_\eps'\|_{L^\infty((0,1))}
  \leq \|f'\|_{L^\infty((0,1))}
$,
but this assumption is not needed for the property
we are proving).
Using that~$g_\eps'$ is bounded for any fixed~$\eps$  
and the estimate $f_\eps \circ g_\eps \leq \eps$,
we thus deduce from~\eqref{bound.form} that
there is a positive constant~$c_{\eps,\delta}$
(again, this constant can be made independent of~$\eps$ 
if $\eps \leq 1$)
such that
$$
  c_{\eps,\delta} \, \|\psi\|_{H_{\gamma_\eps j_\eps/\eps}^1(\Omega)}^2 
  \leq \|\psi\|_{h_\eps}
  \,.
$$
This proves~\eqref{form2}
because $\mathcal{D}$ is dense in $H_{\gamma_\eps j_\eps/\eps}^1(\Omega)$. 
\end{proof}
%

\subsection{The eigenvalue equation}
%
Recall that we denote the eigenvalues of~$T_\eps$ (and hence~$H_\eps$)
by $\mu_n(\Omega_\eps)$ with $n \in \Nat$ ($=\{0,1,\dots\}$). 
The~$(n+1)^\mathrm{th}$ eigenvalue can be characterised 
by the Rayleigh-Ritz variational formula
\begin{equation}\label{minimax.as}
  \mu_n(\Omega_\eps) = 
  \inf_{\stackrel{\dim \mathfrak{L}_n=n+1}
  {\mathfrak{L}_n \subset \Dom(h_\eps)}} 
  \sup_{\psi \in \mathfrak{L}_n} 
  \frac{h_\eps[\psi]}{\|\psi\|^2}
  \,.
\end{equation}
\begin{proposition}\label{Prop.bound}
For any~$n \in \Nat$, 
there exists a positive constant~$C_n$ such that
for all $\eps \leq 1$,
$$
  \mu_n(\Omega_\eps) \leq C_n
  \,.
$$
\end{proposition}
\begin{proof}
Assuming $\eps \leq 1$, we have the following two-sided 
$\eps$- and $t$-independent bound
\begin{equation}\label{2-bound}
  \gamma_-(s)\leq \gamma_\eps(s,t) \leq \gamma_+(s)
\end{equation}
valid for every $(s,t) \in \Omega_+$ with
\begin{equation*}
  \gamma_-(s) :=
  \exp\left(-\frac{(|x_0|+s+1)^2
  +(|y_0|+1)^2}{2}\right)
  , \qquad
  \gamma_+(s) :=
  \exp\left(-\frac{(-|x_0|+s)^2-2|y_0|}{2}\right)
  .
\end{equation*}
Using in addition that 
$g_\eps'=1$ and $f_\eps\circ g_\eps = \eps$ in~$\Omega_+$, 
we obviously have 
$$
  \forall \psi \in C_0^\infty((0,+\infty)) \otimes \{1\}
  \,, \qquad
  \frac{h_\eps[\psi]}{\|\psi\|^2}
  \leq 
  \frac{\displaystyle
  \int_{\Omega_+} 
  (\partial_s \psi)^2 \, \gamma_+(s) \, ds \, dt}
  {\displaystyle \int_{\Omega_+} \psi^2 \, \gamma_-(s) \, ds \, dt}
  \,.
$$
It then follows from~\eqref{minimax.as} that 
the inequality of the proposition holds with the numbers
$$
  C_n := 
  \inf_{\stackrel{\dim \mathfrak{L}_n=n+1}
  {\mathfrak{L}_n \subset C_0^\infty((0,+\infty))}} 
  \sup_{\psi \in \mathfrak{L}_n} 
  \frac{\displaystyle \int_{0}^{+\infty} \psi'(s)^2 \, \gamma_+(s) \, ds}
  {\displaystyle \int_{0}^{+\infty} \psi(s)^2 \, \gamma_-(s) \, ds}
  \,,
$$
which are actually eigenvalues of the one-dimensional operator
$-\gamma_-^{-1} \partial_s \gamma_+ \partial_s$ in 
$\sii_{\gamma_-}((0,+\infty))$,
subject to Dirichlet boundary conditions.
\end{proof}

Let us now fix $n \in \Nat$ and abbreviate the~$(n+1)^\mathrm{th}$
eigenvalue of~$H_\eps$ by $\mu_\eps := \mu_n(\Omega_\eps)$. 
We denote an eigenfunction corresponding to~$\mu_\eps$ by~$\psi_\eps$
and normalise it to~$1$ in $\sii_{\gamma_\eps j_\eps/\eps}(\Omega)$,
\ie, 
\begin{equation}\label{norm}
  \|\psi_\eps\|=1
  \,. 
\end{equation}
for every admissible $\eps > 0$.

The weak formulation of the eigenvalue equation 
$H_\eps \psi_\eps = \mu_\eps \psi_\eps$ reads
\begin{equation}\label{weak}
  \forall \phi \in \Dom(h_\eps) 
  \,, \qquad
  h_\eps(\phi,\psi_\eps) = \mu_\eps \, (\phi,\psi_\eps)
  \,,
\end{equation}
where $(\cdot,\cdot)$ stands for the inner product in 
$\sii_{\gamma_\eps j_\eps/\eps}(\Omega)$
{and $h_\eps(\cdot,\cdot)$ denotes the sesquilinear form
corresponding to~$h_\eps[\cdot]$}, that is $ \forall \phi \in \Dom(h_\eps) $

\begin{align}
   \int_\Omega 
  \left[
  \left(
  \frac{\partial_s \psi_{\eps}}{g_\eps'} 
  - \frac{f_\eps'\circ g_\eps}{f_\eps \circ g_\eps} \, t \, \partial_t \psi_{\eps}
  \right)\left(
  \frac{\partial_s \phi}{g_\eps'} 
  - \frac{f_\eps'\circ g_\eps}{f_\eps \circ g_\eps} \, t \, \partial_t \phi
  \right)
  + 
  \frac{(\partial_t \psi_{\eps})}{(f_\eps\circ g_\eps) }  \frac{(\partial_t \phi)}{(f_\eps\circ g_\eps) }
  \right]
  \gamma_\eps \, g_\eps' \, \frac{f_\eps \circ g_\eps}{\eps} \, ds \, dt
   \qquad
  \label{weform}
  \\
  \notag
  = \mu_\eps \int_\Omega \psi_{\eps}\, \phi \, \gamma_\eps \, g_\eps' \, \frac{f_\eps \circ g_\eps}{\eps} \, ds \, dt
  \,.
\end{align}

\subsection{What happens in \texorpdfstring{$\Omega_+$}{Omega+}}
%
Using $|t| \leq 1$,  we easily verify
\begin{equation}\label{2-bound.bis}
  \forall (s,t) \in \Omega_+ 
  \,, \qquad
  \gamma_\eps(s,t) \geq \rho_\eps(s) \gamma_0(s) 
  \,, 
\end{equation}
where the function
$$
  \rho_\eps(s) := \exp\left(
  -\frac{a_\eps^2+2|x_0| a_\eps + \eps^2+2|y_0|\eps}{2}
  \right)
  \exp(-a_\eps s)
$$
is converging pointwise to~$1$ as $\eps \to 0$.

Choosing $\phi = \psi_\eps$ as a test function in~\eqref{weak}
and using~\eqref{2-bound.bis} 
together with Proposition~\ref{Prop.bound} and~\eqref{norm},
we obtain
\begin{equation}\label{bound.plus}
  \int_{\Omega_+} (\partial_s\psi_\eps)^2 \, \rho_\eps \, \gamma_0 \, ds \, dt
  + \int_{\Omega_+} \frac{(\partial_t\psi_\eps)^2}{\eps^2} 
  \, \rho_\eps \, \gamma_0 \, ds \, dt
  \leq h_\eps[\psi_\eps] = \mu_\eps \|\psi_\eps\|^2 \leq C
  \,.
\end{equation}
Here and in the sequel, we denote by~$C$ a generic
constant which is independent of~$\eps$
and may change its value from line to line.
Writing 
\begin{equation}\label{dec.plus}
  \psi_\eps(s) = \varphi_\eps(s) + \eta_\eps(s,t)
  \,,
\end{equation}
where 
\begin{equation}\label{orthogonal}
  \int_0^1 \eta_\eps(s,t) \, dt = 0
  \qquad \mbox{for a.e.\ } 
  s \in (0,+\infty)
  \,,
\end{equation}
we deduce from the second term on the left hand side of~\eqref{bound.plus}
\begin{equation}\label{eta.plus}
  \pi^2
  \int_{\Omega_+} \eta_\eps^2 \, \rho_\eps\, \gamma_0 \, ds \, dt
  \leq
  \int_{\Omega_+} (\partial_t\eta_\eps)^2 \, \rho_\eps \, \gamma_0 \, ds \, dt
  \leq C \eps^2
  \,.
\end{equation}
Differentiating~\eqref{orthogonal} with respect to~$s$,
we may write
$$
  \int_{\Omega_+} (\partial_s\psi_\eps)^2 \, \rho_\eps \, \gamma_0 \, ds \, dt
  = \int_{\Omega_+} {\varphi_\eps'}^2 \, \rho_\eps \,\gamma_0 \, ds \, dt
  + \int_{\Omega_+} (\partial_s\eta_\eps)^2 \, \rho_\eps \, \gamma_0 \, ds \, dt
$$
and putting this decomposition into~\eqref{bound.plus},
we get from the first term on the left hand side
\begin{equation}\label{varphi'.plus}
  \int_0^{+\infty} {\varphi_\eps'}^2 \, \rho_\eps \, \gamma_0 \, ds 
  \leq C
  \,, \qquad
  \int_0^{+\infty} (\partial_s\eta_\eps)^2 \, \rho_\eps \, \gamma_0 \, ds 
  \leq C
  \,.
\end{equation}

At the same time, from~\eqref{norm} using~\eqref{2-bound.bis},
we obtain  	
\begin{equation}\label{bound.plus.norm}
  \int_{\Omega_+} \varphi_\eps^2 \, \rho_\eps \, \gamma_0 \, ds \, dt
  +\int_{\Omega_+} \eta_\eps^2 \, \rho_\eps \, \gamma_0 \, ds \, dt
  =\int_{\Omega_+} \psi_\eps^2 \, \rho_\eps \, \gamma_0 \, ds \, dt
  \leq \|\psi_\eps\|^2 = 1
  \,,
\end{equation}
where the first equality employs~\eqref{orthogonal}.
Consequently,
\begin{equation}\label{varphi.plus}
  \int_0^{+\infty} \varphi_\eps^2 \, \rho_\eps \, \gamma_0 \, ds 
  \leq 1
  \,.
\end{equation}

Finally, employing the first inequality from~\eqref{varphi'.plus} 
and~\eqref{varphi.plus},
we get
\begin{equation}\label{varphi'.plus.bis}
  \int_0^{+\infty} {(\sqrt{\rho_\eps}\varphi_\eps)'}^2 \, \gamma_0 \, ds 
  \leq C
  \,.
\end{equation}

From~\eqref{varphi.plus} and~\eqref{varphi'.plus.bis},
we see that $\{\sqrt{\rho_\eps}\varphi_\eps\}_{\eps>0}$ 
is a bounded family in $H^1_{\gamma_0}((0,+\infty))$ 
and therefore precompact in the weak topology of this space.
Let~$\varphi_0$ be a weak limit point, 
\ie\ for a decreasing sequence of positive numbers 
$\{\eps_i\}_{i \in \Nat}$
such that $\eps_i \to 0$ as $i \to +\infty$,
\begin{equation}\label{weak.limit}
  \sqrt{\rho_{\eps_i}} \varphi_{\eps_i} 
  \xrightarrow[i\to+\infty]{w} \varphi_0
  \qquad \mbox{in} \qquad H^1_{\gamma_0}((0,+\infty))
  \,.
\end{equation}
Since $H^1_{\gamma_0}((0,+\infty))$ is compactly embedded 
in $\sii_{\gamma_0}((0,+\infty))$, we may assume
\begin{equation}\label{strong.limit}
  \sqrt{\rho_{\eps_i}}\varphi_{\eps_i} \xrightarrow[i\to+\infty]{s} \varphi_0
  \qquad \mbox{in} \qquad \sii_{\gamma_0}((0,+\infty))
  \,.
\end{equation}
%

\subsection{What happens in \texorpdfstring{$\Omega_-$}{Omega-}}
%
Here~$\gamma_\eps$ can be estimated from below just by 
an $\eps$-independent positive number, \eg,
\begin{equation}\label{2-bound.minus}
  \forall (s,t) \in \Omega_- 
  \,, \qquad
  \gamma_\eps(s,t) \geq 
  \exp\left(
  -\frac{(|x_0| +1)^2 + |y_0| +1)^2}{2}
  \right)
  \,.
\end{equation}
On the other hand, we need a lower bound to~$f_\eps$. 
Employing that~$f$ is concave and non-decreasing, 
we can use
\begin{equation}\label{f.bound}
  \forall s \in (-1,0) 
  \,, \qquad
  f_\eps(g_\eps(s)) \geq \eps \, (s+1)
  \,.
\end{equation}
Recall also that $g_\eps'=a_\eps$ on $(-1,0)$.

Choosing $\phi = \psi_\eps$ as a test function in~\eqref{weak}
and using~\eqref{2-bound.minus} and~\eqref{f.bound},
we obtain
\begin{equation}\label{bound.minus}
  \int_{\Omega_-}
  \left(
  \frac{\partial_s\psi_\eps}{a_\eps} 
  - \frac{f_\eps'\circ g_\eps}{f_\eps \circ g_\eps} \, t \, \partial_t\psi_\eps
  \right)^2 
  a_\eps \, (s+1) \, ds \, dt
  + \int_{\Omega_-} \frac{(\partial_t\psi_\eps)^2}{(f_\eps \circ g_\eps)^2} 
  \, a_\eps \, (s+1) \, ds \, dt
  \leq C
  \,.
\end{equation}
Assume~\eqref{Ass}.
Using elementary estimates as in the proof of Proposition~\ref{Prop.form}, 
this inequality implies
\begin{equation}\label{bound.minus.bis}
  \delta \int_{\Omega_-}
  \left(
  \frac{\partial_s\psi_\eps}{a_\eps} 
  \right)^2 a_\eps \, (s+1) \, ds \, dt
  + \left(
  1 - \frac{\delta}{1-\delta} \, \|f_\eps'\|_\infty^2
  \right) 
  \int_{\Omega_-} \frac{(\partial_t\psi_\eps)^2}{(f_\eps \circ g_\eps)^2}
  \, a_\eps \, (s+1) \, ds \, dt
  \leq C
\end{equation}
with any $\delta \in (0,1)$.
We can choose~$\delta$ (independent of~$\eps$ due to~\eqref{Ass.small})
so small that the left hand 
side of~\eqref{bound.minus.bis} is composed of a sum of
two non-negative terms. 
Using in addition $f_\eps \circ g_\eps \leq \eps$,
we thus deduce from~\eqref{bound.minus.bis}
\begin{equation*}
  \frac{1}{a_\eps}
  \int_{\Omega_-}
  (\partial_s\psi_\eps)^2 \, (s+1) \, ds \, dt
  + \frac{a_\eps}{\eps^2}
  \int_{\Omega_-} (\partial_t\psi_\eps)^2
  \, (s+1) \, ds \, dt
  \leq C
  \,.
\end{equation*}
Moreover, it follows from~\eqref{Ass} 
and the convexity bound 
\begin{equation}\label{convexity}
  \forall s \geq 0 
  \,, \qquad
  f(s) \leq f'(0) s
\end{equation}
that
\begin{equation}\label{convexity.corol}
  \eps \leq f'(0) \, a_\eps 
  \,.
\end{equation}
Hence
\begin{equation}\label{bound.minus.bis2}
  \int_{\Omega_-}
  |{\nabla }\psi_\eps|^2 \, (s+1) \, ds \, dt
  \leq C a_\eps
  \,.
\end{equation}

Now we write ($\varphi_\eps$ is constant!)
\begin{equation}\label{dec.minus}
  \psi_\eps(s,t) = \varphi_\eps + \eta_\eps(s,t) 
  \,,
\end{equation}
where 
\begin{equation}\label{orthogonal.minus}
  \int_{\Omega_-} \eta_\eps(s,t) \, (s+1) \, ds\, dt = 0
  \,.
\end{equation}
Then we deduce from~\eqref{bound.minus.bis2}
\begin{equation}\label{eta.minus}
  \pi^2
  \int_{\Omega_-} \eta_\eps^2 \, (s+1) \, ds \, dt
  \leq
  \int_{\Omega_-} |{\nabla }\eta_\eps|^2 \, (s+1) \, ds \, dt
  \leq C a_\eps
  \,.
\end{equation}
Note that~$\pi^2$ is indeed the minimum between the first
non-zero Neumann eigenvalue in the interval of unit length
and the first non-zero Neumann eigenvalue in the unit disk.

At the same time, from~\eqref{norm} using~\eqref{2-bound.minus}
and~\eqref{f.bound}, we obtain  	
\begin{equation}\label{bound.minus.norm}
  \int_{\Omega_-} \varphi_\eps^2 \, a_\eps \, (s+1) \, ds \, dt
  +\int_{\Omega_-} \eta_\eps^2 \, a_\eps \, (s+1) \, ds \, dt
  = \int_{\Omega_-} \psi_\eps^2 \, a_\eps \, (s+1) \, ds \, dt
  \leq C
  \,,
\end{equation}
where the first equality employs~\eqref{orthogonal.minus}. 
Consequently, recalling that~$\varphi_\eps$ is constant,
\begin{equation}\label{varphi.minus}
  \varphi_\eps^2 \, a_\eps  
  \leq C
  \qquad \mbox{on} \quad \Omega_-
  \,.
\end{equation}
%

\subsection{The limiting eigenvalue equation in \texorpdfstring{$\Omega_+$}{Omega+}}
%
Now we consider~\eqref{weak} 
for the sequence $\{\eps_i\}_{i\in\Nat}$ 
and a test function $\phi(s,t) = \varphi(s)$,
where $\varphi \in C_0^\infty(\Real)$
is such that $\varphi'=0$ on $[-1,0]$,
and take the limit $i \to +\infty$.

We shall need a lower bound analogous 
to the upper bound~\eqref{convexity.corol}.
From the fundamental theorem of calculus,
we deduce 
\begin{equation}\label{Taylor}
  \forall s \in [0,a_\eps] \,, \qquad
  f(s) \geq \big( \essinf_{(0,a_\eps)} f' \big) \, s
  \,.
\end{equation}
Note that the infimum cannot be zero unless~$f$ is trivial
(we assume from the beginning $\eps < \sup f$
and that~$f$ is non-decreasing)
and that it converges to $f'(0)>0$ as $\eps \to 0$. 
Consequently, for all sufficiently small~$\eps$, we have
\begin{equation}\label{convexity.corol.lower}
  \eps \geq \frac{1}{2} \, f'(0) \, a_\eps 
  \,.
\end{equation}

At the same time, in analogy with~\eqref{2-bound.bis},
we have
\begin{equation}\label{2-bound.new}
  \forall (s,t) \in \Omega_+ 
  \,, \qquad
  \gamma_\eps(s,t) \leq c_\eps \rho_\eps(s) \gamma_0(s) 
  \,, 
\end{equation}
where
$$
  c_\eps := \exp\left(
  \frac{2 |x_0| a_\eps + \eps^2 + 2|y_0|\eps}{2}
  \right)
$$
is converging to~$1$ as $\eps \to 0$.

We first look at the right hand side of~\eqref{weak}.
Using the decompositions~\eqref{dec.plus} and~\eqref{dec.minus},
we have
\begin{multline*}
  (\varphi,\psi_{\eps})
  = \int_{\Omega_-} \varphi \, \varphi_{\eps} \,
  \gamma_\eps \, a_\eps \, \frac{f_\eps \circ g_\eps}{\eps} \, ds \, dt
  + \int_{\Omega_-} \varphi \, \eta_\eps \, 
  \gamma_\eps \, a_\eps \, \frac{f_\eps \circ g_\eps}{\eps} \, ds \, dt
  \\
  + \int_{\Omega_+} \varphi \, \varphi_\eps \,
  \gamma_\eps \, ds \, dt
  + \int_{\Omega_+} \varphi \, \eta_\eps \, 
  \gamma_\eps \, ds \, dt
  \,.
\end{multline*}
Estimating $\gamma_\eps \leq 1$ and using~\eqref{convexity}
and~\eqref{convexity.corol.lower},
we get
$$
  \left|
  \int_{\Omega_-} \varphi \, \eta_\eps \, 
  \gamma_\eps \, a_\eps \, \frac{f_\eps \circ g_\eps}{\eps} \, ds \, dt
  \right|
  \leq 
  \frac{a_\eps^2}{\eps} f'(0)
  \int_{\Omega_-} |\varphi| \, |\eta_\eps| \, 
  (s+1) \, ds \, dt
  \leq 
  2 a_\eps 
  \int_{\Omega_-} |\varphi| \, |\eta_\eps| \, 
  (s+1) \, ds \, dt
  \,,
$$ 
where the right hand side tends to zero as $\eps \to 0$
due to the Schwarz inequality and~\eqref{eta.minus}. 
At the same time, 
recalling that~$\varphi_\eps$ is constant in~$\Omega_-$,
$$
  \left|
  \int_{\Omega_-} \varphi \, \varphi_\eps \, 
  \gamma_\eps \, a_\eps \, \frac{f_\eps \circ g_\eps}{\eps} \, ds \, dt
  \right|
  \leq 
  \frac{a_\eps^2}{\eps} f'(0)
  \int_{\Omega_-} |\varphi| \, |\varphi_\eps| \, 
  (s+1) \, ds \, dt
  \leq 
  2 a_\eps |\varphi_\eps|
  \int_{\Omega_-} |\varphi| \,  
  (s+1) \, ds \, dt
  \,,
$$ 
where the right hand side tends to zero as $\eps \to 0$
due to~\eqref{varphi.minus}.
Using~\eqref{2-bound.new}, we also get 
$$
  \left|
  \int_{\Omega_+} \varphi \, \eta_\eps \, 
  \gamma_\eps \, ds \, dt
  \right|
  \leq 
  c_\eps
  \int_{\Omega_+} |\varphi| \, |\eta_\eps| \, 
  \rho_\eps \, \gamma_0 \, ds \, dt
  \leq
  c_\eps
  \sqrt{\int_{\Omega_+} \varphi^2 \, 
  \gamma_0 \, ds \, dt}
  \sqrt{\int_{\Omega_+} \eta_\eps^2 \, 
  \rho_\eps \, \gamma_0 \, ds \, dt}
  \,,
$$ 
where the right hand side tends to zero as $\eps \to 0$
due to~\eqref{eta.plus}. 
Finally, we write 
$$
  \int_{\Omega_+} \varphi \, \varphi_{\eps_i} 
  \gamma_{\eps_i} ds \, dt
  = \int_{\Omega_+} \varphi \, \varphi_{\eps_i} 
  \sqrt{\rho_{\eps_i}} \gamma_{0} \, ds \, dt
  + \int_{\Omega_+} \varphi \, \varphi_{\eps_i} 
  \sqrt{\rho_{\eps_i}} \gamma_{0} \,
  \left(\frac{\gamma_{\eps_i}}{\sqrt{\rho_{\eps_i}} \gamma_0}
  -1\right)
  \, ds \, dt
  \,.
$$
Here the first term on the right hand side converges to 
$
  \int_{\Omega_+} \varphi \, \varphi_{0} 
  \, \gamma_{0} \, ds \, dt
$
as $i \to +\infty$ due to~\eqref{weak.limit},
while the second term vanishes in the limit because of 
\begin{multline*}
  \left|
  \int_{\Omega_+} \varphi \, \varphi_{\eps_i} 
  \sqrt{\rho_{\eps_i}} \gamma_{0} \,
  \left(\frac{\gamma_{\eps_i}}{\sqrt{\rho_{\eps_i}} \gamma_0}
  -1\right)
  \, ds \, dt
  \right|
  \\
  \leq
  \sqrt{
  \int_{\Omega_+} \varphi^2 
  \, \gamma_{0} \,
  \left(\frac{\gamma_{\eps_i}}{\sqrt{\rho_{\eps_i}} \gamma_0}
  -1\right)^2
  \, ds \, dt
  }
  \sqrt{
  \int_{\Omega_+} \varphi_{\eps_i}^2 
  \rho_{\eps_i} \gamma_{0} \,
  \, ds \, dt
  }
  \,.
\end{multline*}
Indeed the second term on the right hand side is bounded 
by~\eqref{bound.plus.norm}, 
while first term tends to zero as $i \to +\infty$
by the dominated convergence theorem.
Summing up,
\begin{equation}\label{lim.right}
  \lim_{i \to +\infty} (\varphi,\psi_{\eps_i})
  = \int_{0}^{+\infty} \varphi \, \varphi_{0} 
  \, \gamma_{0} \, ds 
  \,.
\end{equation}

Employing that the test function~$\varphi$ is constant on $[-1,0]$
and the decomposition~\eqref{dec.plus},
we have
$$
  h_\eps(\varphi,\psi_{\eps})
  = \int_{\Omega_+} \varphi' \, \varphi_\eps' \,
  \gamma_\eps \, ds \, dt
  + \int_{\Omega_+} \varphi' \, \partial_s\eta_\eps \, 
  \gamma_\eps \, ds \, dt
  \,.
$$  
Here the first term on the right hand side can treated 
in the same way as above with the conclusion
$$
  \int_{\Omega_+} \varphi' \, \varphi_{\eps_i}' \,
  \gamma_{\eps_i} \, ds \, dt
  \xrightarrow[i \to +\infty]{}
  \int_{\Omega_+} \varphi' \, \varphi_{0}' \,
  \gamma_{0} \, ds \, dt
  = \int_{0}^{+\infty} \varphi' \, \varphi_{0}' \,
  \gamma_{0} \, ds 
  \,,
$$
while we integrate by parts to handle the second term,
$$
  \int_{\Omega_+} \varphi' \, \partial_s\eta_\eps \, 
  \gamma_\eps \, ds \, dt
  = - \int_{\Omega_+} \varphi'' \, \eta_\eps \, 
  \gamma_\eps \, ds \, dt
  - \int_{\Omega_+} \varphi' \, \eta_\eps \, 
  \partial_s\gamma_\eps \, ds \, dt
  \,.
$$
Notice that the boundary terms vanish because~$\varphi$
has a compact support in~$\Real$ and $\varphi'(0)=0$. 
As above, the first term on the right hand side
vanishes as $\eps \to 0$ due to~\eqref{eta.plus}. 
Similarly, 
$$
\begin{aligned}
  \left|
  \int_{\Omega_+} \varphi' \, \eta_\eps \, 
  \partial_s\gamma_\eps \, ds \, dt
  \right|
  &\leq 
  c_\eps
  \int_{\Omega_+} |\varphi| \, |\eta_\eps| \, 
  \rho_\eps \, \gamma_0 \, (x_0+s+a_\eps)
  \, ds \, dt
  \\
  &\leq
  c_\eps
  \sqrt{\int_{\Omega_+} \varphi^2 \, 
  \gamma_0 \, (x_0+s+a_\eps)^2 \, ds \, dt}
  \sqrt{\int_{\Omega_+} \eta_\eps^2 \, 
  \rho_\eps \, \gamma_0 \, ds \, dt}
  \,,
\end{aligned}
$$ 
where the right hand side tends to zero as $\eps \to 0$
due to~\eqref{eta.plus}. 
Summing up,
\begin{equation}\label{lim.left}
  \lim_{i \to +\infty} h_{\eps_i}(\varphi,\psi_{\eps_i})
  = \int_{0}^{+\infty} \varphi' \, \varphi_{0}' 
  \, \gamma_{0} \, ds 
  \,.
\end{equation}

Since the set of functions $\varphi \in C_0^\infty(\Real)$ 
satisfying $\varphi'(0)=0$ is a core for the form 
domain of the operator~$T_0$, 
we conclude from~\eqref{lim.left} and~\eqref{lim.right}
that~$\varphi_0$ belongs to~$\Dom(T_0)$
and solves the one-dimensional problems
\begin{equation}\label{1Dproblems}
\begin{aligned}
  T_0\varphi_0 &= \mu_0^+ \varphi_0
  \,, \qquad &
  \mu_0^+ &:= \limsup_{i\to+\infty} \mu_{\eps_i}
  \,, 
  \\
  T_0\varphi_0 &= \mu_0^- \varphi_0
  \,, \qquad &
  \mu_0^- &:= \liminf_{i\to+\infty} \mu_{\eps_i}
  \,.
\end{aligned}
\end{equation}
If $\varphi_0 \not=0$ on $(0,+\infty)$, 
then~$\mu_0^\pm$ must coincide
with some eigenvalues of~$T_0$.  
It remains to check that indeed $\varphi_0 \not=0$ on $(0,+\infty)$.

\subsection{The limiting problem in \texorpdfstring{$\Omega_-$}{Omega-}: 
a crucial step}
%
Define
$$
  \Omega_-' := (-1/2,0)\times(0,1)
  \,, \qquad 
  \Omega_+' := (0,1/2)\times(0,1)
  \,, \qquad
  \Omega' := (-1/2,1/2)\times(0,1)
  \,.
$$
From~\eqref{bound.minus.bis2} and~\eqref{bound.minus.norm},
we respectively have 
\begin{equation}
  \int_{\Omega_-'} |{\nabla }\psi_\eps|^2 \, ds \, dt \leq 2C a_\eps
  \,, \qquad
  \int_{\Omega_-'} \psi_\eps^2 \, ds \, dt \leq \frac{2C}{a_\eps} 
  \,.
\end{equation}
At the same time, denoting 
$
  m_0 := \min_{[0,1/2]} \gamma_0
$ 
and assuming $\eps \leq 1$,
from~\eqref{bound.plus} and~\eqref{bound.plus.norm},
we respectively get 
\begin{equation}
  \int_{\Omega_+'} |{\nabla }\psi_\eps|^2 \, ds \, dt 
  \leq \frac{C}{m_0 \, \rho_\eps(1/2)}
  \,, \qquad
  \int_{\Omega_+'} \psi_\eps^2 \, ds \, dt 
  \leq \frac{1}{m_0 \, \rho_\eps(1/2)} 
  \,.
\end{equation}
Consequently, $\psi_\eps \in H^1(\Omega')$ for any $\eps \leq 1$
(although, in principle, 
$\|\psi_\eps\|_{H^1(\Omega')}$ might not be uniformly bounded in~$\eps$).

It follows that the boundary values $\psi_\eps(0-,t)$ and $\psi_\eps(0+,t)$
exist in the sense of traces in~$\Omega_-'$ and~$\Omega_+'$, respectively,
and they must be equal as functions of~$t$ in $\sii((0,1))$.  
Using the decompositions~\eqref{dec.plus} and~\eqref{dec.minus},
we therefore have, for almost every $t \in (0,1)$, 
\begin{multline*}
  [\varphi_\eps(0-) - \varphi_\eps(0+)]^2
  = \left[
  \eta_\eps(0+,t) - \eta_\eps(0-,t) 
  \right]^2
  \leq 2\,[\eta_\eps(0+,t)]^2 + 2\,[\eta_\eps(0-,t)]^2 
  \\
  \leq 2\,C \int_0^{1/2} 
  \left([\eta_\eps(s,t)]^2 + [\partial_s\eta_\eps(s,t)]^2\right) ds 
  +  2\,C \int_{-1/2}^0 
  \left([\eta_\eps(s,t)]^2 + [\partial_s\eta_\eps(s,t)]^2\right) ds 
  \,,
\end{multline*}
where~$C$ is a constant coming from the Sobolev embedding theorem.
Recall that~$\varphi_\eps$ is constant on~$(-1,0)$
and $\varphi_\eps \in H^1((0,1/2)) \hookrightarrow C^0([0,1/2])$;
more specifically, the first inequality of~\eqref{varphi'.plus} 
and~\eqref{varphi.plus} respectively yield
\begin{equation}\label{varphi.H1}
  \int_0^{1/2} {\varphi_\eps'}^2 \, ds 
  \leq \frac{C}{m_0 \, \rho_\eps(1/2)}
  \,, \qquad
  \int_0^{1/2} \varphi_\eps^2 \, ds 
  \leq \frac{1}{m_0 \, \rho_\eps(1/2)}
  \,.
\end{equation}
Integrating with respect to~$t$ above, we deduce 
$$
  [\varphi_\eps(0-) - \varphi_\eps(0+)]^2
  \leq 2\,C \int_{\Omega_+'} 
  \left[\eta_\eps^2 + (\partial_s\eta_\eps)^2\right] ds \, dt 
  + 2\,C \int_{\Omega_-'} 
  \left[\eta_\eps^2 + (\partial_s\eta_\eps)^2\right] ds \, dt 
  \,.
$$
Applying~\eqref{eta.plus}, 
the second inequality of~\eqref{varphi'.plus} and~\eqref{eta.minus},
we may write 
\begin{equation}\label{trace.trick}
  [\varphi_\eps(0-) - \varphi_\eps(0+)]^2 \leq C
  \,,
\end{equation}
where~$C$ is a constant (different from the above) independent of~$\eps$, 
provided that~\eqref{Ass.small} holds.
Finally, applying~\eqref{varphi.H1} and the Sobolev embedding
$H^1((0,1/2)) \hookrightarrow C^0([0,1/2])$,
we deduce from~\eqref{trace.trick} the following improvement
upon~\eqref{varphi.minus}
\begin{equation}\label{varphi.minus.improve}
  \varphi_\eps^2 
  \leq C
  \qquad \mbox{on $\Omega_-$}
  \,.
\end{equation}
%

\subsection{As \texorpdfstring{$\eps \to 0$}{epsilonto0} only 
\texorpdfstring{$\Omega_+$}{Omega+} matters: 
convergence of eigenvalues and eigenfunctions}
%
Estimate \eqref{varphi.minus.improve}~provides a crucial information 
whose significance consists in that
what happens in~$\Omega_-$ is insignificant.
\begin{proposition}\label{Prop.norm}
One has
$$  
  \|\psi_{\eps_i}\| \xrightarrow[i \to +\infty]{} 
  \|\varphi_0\|_{\sii_{\gamma_0}((0,+\infty))}
  \,.
$$
\end{proposition}
\begin{proof}
We have 
\begin{align*}
  \|\psi_{\eps}\|^2
  = & \int_{\Omega_-} \varphi_{\eps}^2 \,
  \gamma_\eps \, a_\eps \, \frac{f_\eps \circ g_\eps}{\eps} \, ds \, dt
  + \int_{\Omega_-} \eta_\eps^2 \, 
  \gamma_\eps \, a_\eps \, \frac{f_\eps \circ g_\eps}{\eps} \, ds \, dt
  \\
  & + \int_{\Omega_-} 2 \, \varphi_{\eps} \, \eta_\eps \,
  \gamma_\eps \, a_\eps \, \frac{f_\eps \circ g_\eps}{\eps} \, ds \, dt
  \\
  & + \int_{\Omega_+} \varphi_{\eps}^2 \,
  \gamma_\eps \, ds \, dt
  + \int_{\Omega_+} \eta_\eps^2 \, 
  \gamma_\eps \, ds \, dt
  + \int_{\Omega_+} 2 \, \varphi_{\eps} \, \eta_\eps \,
  \gamma_\eps \, ds \, dt
  \,.
\end{align*}
The right hand side of the first line 
together with the mixed term on the second line
goes to zero as $\eps \to 0$.
Indeed, recalling~\eqref{convexity}, \eqref{convexity.corol.lower}
and $\gamma_\eps \leq 1$,
$$
  \int_{\Omega_-} \varphi_{\eps}^2 \,
  \gamma_\eps \, a_\eps \, \frac{f_\eps \circ g_\eps}{\eps} \, ds \, dt
  \leq 
  2 \, a_\eps \, \varphi_{\eps}^2 \int_{\Omega_-} (s+1) \, ds \, dt
  \xrightarrow[\eps \to 0]{} 0
$$
due to~\eqref{varphi.minus.improve};
$$
  \int_{\Omega_-} \eta_{\eps}^2 \,
  \gamma_\eps \, a_\eps \, \frac{f_\eps \circ g_\eps}{\eps} \, ds \, dt
  \leq 
  2 \, a_\eps \int_{\Omega_-} \eta_{\eps}^2 \, (s+1) \, ds
  \xrightarrow[\eps \to 0]{} 0
$$
due to~\eqref{eta.minus};
and the mixed term goes to zero by the Schwarz inequality.
Similarly, recalling~\eqref{2-bound.new},
$$
  \int_{\Omega_+} \eta_{\eps}^2 \,
  \gamma_\eps \, ds \, dt
  \leq 
  c_\eps \int_{\Omega_+} \eta_{\eps}^2 \, \rho_\eps \gamma_0 \, ds \, dt
  \xrightarrow[\eps \to 0]{} 0
$$
due to~\eqref{eta.plus};
while the Schwarz inequality yields
$$
  \left|
  \int_{\Omega_+} 2\, \varphi_{\eps} \, \eta_\eps \,
  \gamma_\eps \, ds \, dt
  \right|
  \leq 
  2 c_\eps 
  \sqrt{
  \int_{\Omega_+} \eta_{\eps}^2 \, \rho_\eps \gamma_0 \, ds \, dt
  }
  \sqrt{
  \int_{\Omega_+} \varphi_{\eps}^2 \, \rho_\eps \gamma_0 \, ds \, dt
  }
  \xrightarrow[\eps \to 0]{} 0
  \,,
$$
where the second square root is bounded in~$\eps$ 
due to~\eqref{varphi.plus}.
Finally, we write
$$
  \int_{\Omega_+} \varphi_{\eps_i}^2 \,
  \gamma_{\eps_i} \, ds \, dt
  = \int_{\Omega_+} \varphi_{\eps_i}^2 \,
  \rho_{\eps_i} \gamma_0 \, ds \, dt
  + \int_{\Omega_+} \varphi_{\eps_i}^2 \,
  (\gamma_{\eps_i}-\rho_{\eps_i} \gamma_0) \, ds \, dt
$$
and observe that the first term on the right hand side
tends to the desired result
$\|\varphi_0\|_{\sii_{\gamma_0}((0,+\infty))}^2$
as $i \to+ \infty$ by the strong convergence~\eqref{strong.limit},
while the second term vanishes in the limit.
In more detail,
\begin{multline*}
  \left|  
  \int_{\Omega_+} \varphi_{\eps_i}^2 \,
  (\gamma_{\eps_i}-\rho_{\eps_i} \gamma_0) \, ds \, dt
  \right|
  = \left|  
  \int_{\Omega_+} 
  \left(
  \varphi_{\eps_i}^2 \, \rho_{\eps_i}
  - \varphi_{0}^2 + \varphi_{0}^2
  \right)
  \left(\frac{\gamma_{\eps_i}}{\rho_{\eps_i}}-\gamma_0\right) \, ds \, dt
  \right|
  \\
  \leq 
  \int_{\Omega_+} 
  \left|
  \varphi_{\eps_i}^2 \, \rho_{\eps_i}
  - \varphi_{0}^2 
  \right|
  \left(c_{\eps_i}\gamma_0+\gamma_0\right) \, ds \, dt
  + \int_{\Omega_+} 
  \varphi_{0}^2
  \left(\frac{\gamma_{\eps_i}}{\rho_{\eps_i}}-\gamma_0\right) \, ds \, dt
  \,,
\end{multline*}
where the the first term after the inequality tends to zero as $i \to +\infty$
by the strong convergence again, while the second term vanishes
by the dominated convergence theorem.
\end{proof}

It follows from Proposition~\ref{Prop.norm} that $\varphi_0 \not= 0$,
so that it is indeed an eigenfunction of~$T_0$
due to~\eqref{1Dproblems}. 
In particular, $\mu_0^+ = \mu_0^-$.

Now, let~$\hat{\psi}_\eps$ be a normalised eigenfunction
corresponding to possibly another eigenvalue 
$\hat{\mu}_\eps := \mu_m(\eps)$.   
Again, we use the decompositions~\eqref{dec.plus} and~\eqref{dec.minus}
and distinguish the individual components by tilde.
In the same way as we proved Proposition~\ref{Prop.norm},
we can establish
\begin{proposition}\label{Prop.inner}
One has
$$  
  (\psi_{\eps_i},\hat{\psi}_{\hat\eps_j}) 
  \xrightarrow[i,j \to +\infty]{} 
  (\varphi_0,\hat{\varphi}_0)_{\sii_{\gamma_0}((0,+\infty))}
  \,.
$$
\end{proposition}

If $m \not= n$, then $(\psi_{\eps_i},\hat{\psi}_{\hat\eps_j})=0$
and thus $(\varphi_0,\hat{\varphi}_0)_{\sii_{\gamma_0}((0,+\infty))}=0$.
Hence~$\varphi_0$ and~$\hat{\varphi}_0$ correspond to
distinct eigenvalues of~$T_0$. 
In particular, $\varphi_0$~is an eigenfunction 
corresponding to the $(n+1)^\mathrm{th}$ eigenvalue~$\nu_n$ of~$T_0$.
Since we get this result for \emph{any} weak limit 
point of $\{\varphi_\eps\}_{\eps>0}$, 
we have the convergence results actually in $\eps \to 0$
(no need to pass to subsequences).

This completes the proof of Theorem~\ref{Thm.convergence}.

\vskip 1cm
\noindent \textsc{Acknowledgements}. 
This paper was partially supported by the grants PRIN 2012 
``Elliptic and parabolic partial differential equations: 
geometric aspects, related inequalities, and applications'', 
FIRB 2013  ``Geometrical and qualitative aspects of PDE's'', 
and  STAR 2013 ``Sobolev-Poincar\'e inequalities: embedding constants, 
stability issues, nonlinear eigenvalues'' (SInECoSINE).  
{D.K.\ was partially supported
by the project RVO61389005 and the GACR grant No.\ 14-06818S.}

%

\end{document}